\documentclass[12pt]{article}

\usepackage{amsmath}
\usepackage{upgreek}
\usepackage{enumitem} 
\usepackage{amsthm}
\usepackage{graphicx}
\usepackage{xcolor}
\usepackage{esdiff}
\usepackage{mathtools}
\usepackage{amssymb}
\usepackage{ marvosym }
\usepackage{amsmath,amssymb,amsfonts}
\usepackage{algorithmic}
\usepackage{graphicx}
\usepackage{textcomp}
\usepackage{xcolor}
\usepackage{float}
\usepackage[
backend=biber,
style=alphabetic,
sorting=ynt
]{biblatex}
\DeclareUnicodeCharacter{03B8}{HERE!HERE!}

\begin{document}

\title{Data Clustering and Visualization with Recursive Max $k$-Cut Algorithm}
\author{An Ly, Raj Sawhney, Marina Chugunova}
\date{July 23rd, 2024}
% \address{Claremont Graduate University}
% \email{an.ly@cgu.edu, raj.sawhney@cgu.edu, marina.chugunova@cgu.edu}
\maketitle

\begin{abstract}
In this article, we continue our analysis for a novel recursive modification to the Max $k$-Cut algorithm using semidefinite programming as its basis, offering an improved performance in vectorized data clustering tasks. Using a dimension relaxation method, we use a recursion method to enhance density of clustering results. Our methods provide advantages in both computational efficiency and clustering accuracy for grouping datasets into three clusters, substantiated through comprehensive experiments.
\end{abstract}

\section{Introduction}
With the increasing amount of biomedical articles being published every year, researchers have looked into clustering them based on features such as citations, topics, and other measures of similarity as researched in \cite{comparingrelatednessmeasures} or \cite{clusterbycitation}. 
Clustering these documents is crucial for information retrieval and for analysis in not only the biomedical field but also in multiple modern research projects. 
In order to accurately cluster these articles into usable groups, many algorithms have been developed and tested by researchers over the years. 
With over two million biomedical publications in existence, and more being published every year, several similarity based approaches have been used to accurately cluster these articles into groups as discussed in \cite{twomillionbiomedicalpublications} and \cite{contentbasedpublicationalgorithm} by extracting keywords as done in \cite{unsupervisedkeyphraseextraction}.

MaxCut and Max $k$-Cut algorithms have been well-researched and used to cluster vectorized datasets. 
These algorithms have been extensively researched through the use of semidefinite programming, random strategies, adaptive searches, and other methods as seen in \cite{algorithmcomparison}, \cite{randomheuristics}, \cite{karger1998approximate}, and \cite{computationalapproach}.

Other approaches for solving max $k$-cut specifically have included polynomial-time Gaussian sampling-based algorithms by \cite{shinde2021memoryefficient}, branch-and-bound algorithms by \cite{Lu_Deng_2021}, multiple operator heuristics by \cite{Ma_Hao_2016}, and several others as seen in \cite{Hojny_Joormann_Lüthen_Schmidt_2020} and \cite{Zhu_Guo_2011}.

Here, we introduce both a recursive application and a higher-dimension relaxation method to a semidefinite program for the Max $k$-Cut Algorithm created by \cite{Frieze_Jerrum_1997}.
The Max $k$-Cut problem, similar to MaxCut, has been extensively discussed in research such as \cite{Kann_Khanna_Lagergren_Panconesi_1997}, \cite{Arora_Lund_Motwani_Sudan_Szegedy}, \cite{Delorme_Poljak_1993}, and \cite{PAPADIMITRIOU1991425}.
Specifically, we apply this to the case where $k = 3$ which has an expected cut value of at least $80\%$ of the optimal value as discussed in Section~\ref{maxkcutsummary}. 
 
\section{A Summary of Max $k$-Cut}
\label{maxkcutsummary}
For a MaxCut algorithm, given a dataset we would like to split the index set (each index corresponding to a paragraph) to maximize the dissimilarity between the resulting subsets.

In other words, given a matrix of similarity weights $W = \{w_{ij}\}$ where $i, j = 1, \cdots, n$ (where $n$ is the number of datapoints), we want to split the index set $i = 1, \cdots, n$ into $k$ sets $A_1, A_2, \cdots, A_k$ such that $\sum\limits_{i < j} w_{ij}$ is maximized for $l < m$ where $l,m = 1, \cdots, k$ and $i \in A_l, j \in A_m$. 
The expected value of total dissimilarity cut obtained by GWA will be at least $83.6\%$ of the optimal value as proved by Goemans and Williamson in \cite{Goemans_Williamson_2004}. 
In this proof, they utilize a randomized $(\rho - \epsilon)$-approximation algorithm for this estimated value.

Frieze and Jerrum utilize semidefinite programming for Max $k$-Cut in \cite{Frieze_Jerrum_1997} to achieve an expected cut value of at least $80.02\%$.
We include a summarized version of the algorithm below. 

The algorithm begins by creating an equilateral simplex in $R^{k-1}$, the $k$-dimensional space (where $k$ is the number of desired clusters) with vertices $b_1, \cdots, b_k$. 
After doing so, they find the centroid $c_k$ of this simplex and define $a_i = b_i - c_k$ for the values of $1 \leq i \leq k$. 
Then each $a_i$ vector is normalized.  
Next, assign $y_j := a_i$ where $j = 1, \cdots, n$.
The assignment algorithm is described later. 

Note that 
\begin{equation}
    \langle a_i, a_j \rangle = -\frac{1}{k-1}
\end{equation} 
for all $i \neq j$. 
Then we formulate the optimization problem for Max $k$-Cut as: \begin{equation}
    \max \frac{k-1}{k} \sum_{i < j} w_{ij}(1-\langle y_i, y_j \rangle)
\end{equation} subject to the constraints that $y_j \in \{a_1, a_2, \cdots, a_k\}$. 

To obtain the semidefinite programming relaxation, we replace each $y_i$ with a corresponding vector $v_i$ where $v_i \in S_{n-1}$, the $n$-dimensional unit sphere. 
This relaxation results in the following: \begin{equation}
    \max \frac{k-1}{k}\sum_{i < j} w_{ij}(1 - \langle v_i, v_j \rangle)
\end{equation} subject to the constraints that $v_j \in S_{n-1}$ for $j = 1, \cdots, n$ and $\langle v_i, v_j \rangle \geq -\frac{1}{k-1}$ for all $i \neq j$. 

We apply this to our dataset in the following manner:\\
1) Solve the relaxed optimization problem to obtain our optimal vectors $v_i$.\\
2) Choose $k$ random vectors $z_1, \cdots, z_k$ using a multivariate normal distribution.\\
3) Calculate $\langle v_i, z_l \rangle$ for all $l = 1, \cdots, k$.\\
4) Assign $v_i$ to cluster $j$ where $j = \text{argmax}_{l} \langle v_i, z_l \rangle$. 

The analytical proof for why this process achieves an expected cut value of at least $80.02\%$ of the optimal value can be found in \cite{Frieze_Jerrum_1997}. 
Klerk et. al achieved an expected cut value of $83.6\%$ \cite{Klerk_Pasechnik_Warners_2004}. 

For our research, further refinements have been made to our model using additional knowledge found in \cite{Fakhimi_Validi_Hicks_Terlaky_Zuluaga_2023} which researches relaxations on max $k$-cut problems and \cite{newman:OASIcs.SOSA.2018.13} detailing a rounding algorithm for standard semidefinite programming relaxation.

\section{Recursive Generalization of Max $k$-Cut Algorithm}
\label{recursionalgorithm}
We define our recursion algorithm in the following way:\\
1) Apply Frieze and Jerrum's algorithm as discussed in Section~\ref{maxkcutsummary}. \\
2) Repeat Steps 2 through 4 for a total of 200 iterations.\\
3) At each iteration, compute the dissimilarity between clusters and the dissimilarity inside of clusters. Update the optimal partition only if the dissimilarity between clusters is larger and the dissimilarity inside clusters is smaller. Note that the dissimilarity inside clusters is defined as: $\sum_{l<m}\sum_{i<j} w_{ij}$ where $i,j \in A_l$.\\
4) Calculate the 2D-Principal Component Analysis of the $V$ matrix whose columns are the calculated $v_i$ vectors. \\
5) Use the results of Step 4 as the new input of the algorithm.\\
6) Repeat Steps 1 through 5 as desired. 

Utilizing principal component analysis to visualize the obtained matrix $V$, we can see in our experimental results how the $v_i$ vectors begin to cluster, a possible graphical explanation for why the lower bound of our Max 3-cut value is at least $80\%$ of the optimal cut value in comparison to $33\%$ from random assignment algorithms. 
For Step 4, we construct the new dataset by applying PCA to $V$ and using each resulting row as a datapoint.

We apply this algorithm to three datasets: one hand-generated moonset dataset as seen in Section~\ref{recursiveresults}, a collection of spike timestamps of simulated brain waves as discussed in Section~\ref{brainwavedataset}, and a collection of paragraphs from biomedical articles about amodiaquine in Section~\ref{amodiaquinedataset}. 

\subsection{Results of Recursive GWA}
\label{recursiveresults}
Here, we analyze the results of our recursive Max $k$-Cut algorithm on a moonset dataset which contains 100 points in total. 

\begin{figure}
\centering
        \includegraphics[width = 4.3cm]{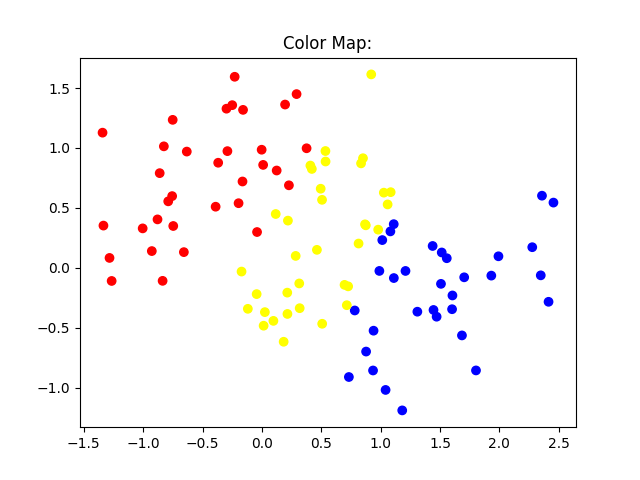}
        \includegraphics[width = 4.3cm]{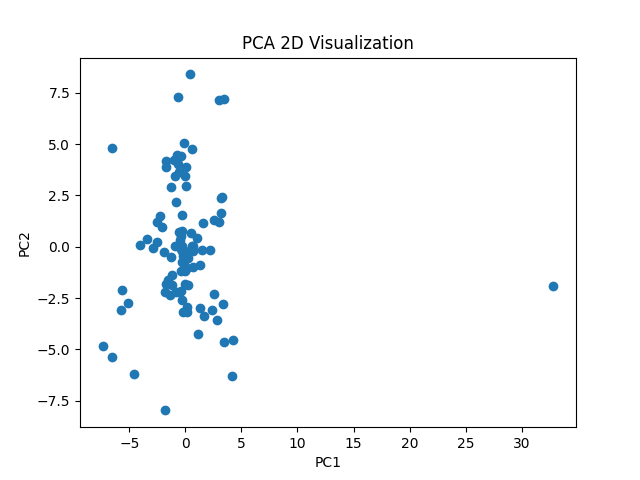}\\
        \includegraphics[width = 4.3cm]{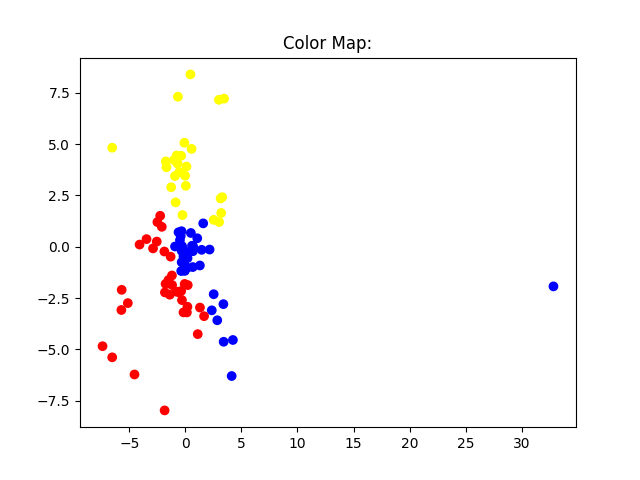}
        \includegraphics[width = 4.3cm]{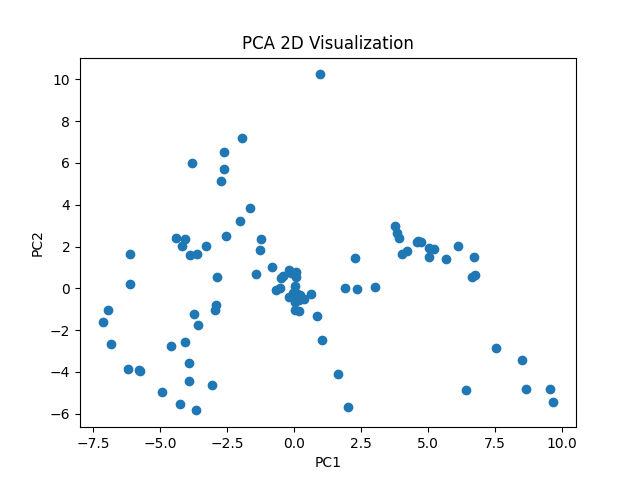}\\
        \includegraphics[width = 4.3cm]{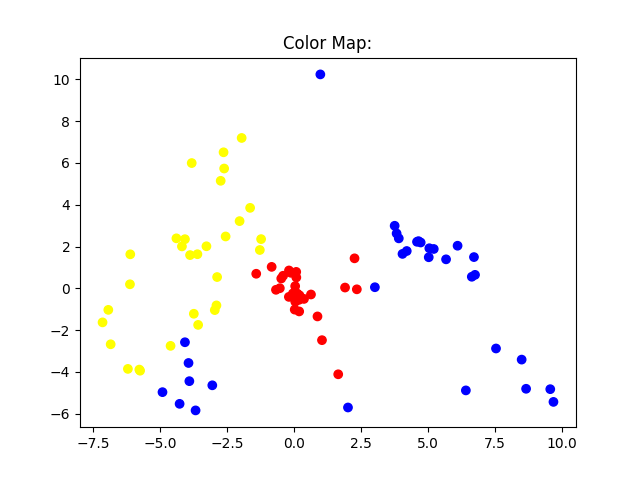}
        \includegraphics[width = 4.3cm]{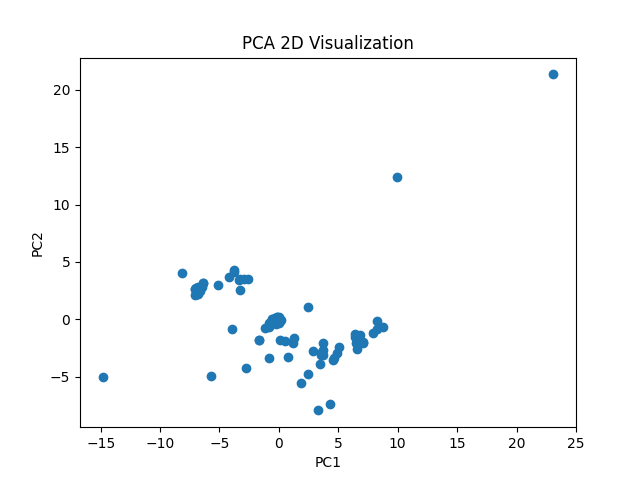}\\
        \includegraphics[width = 4.3cm]{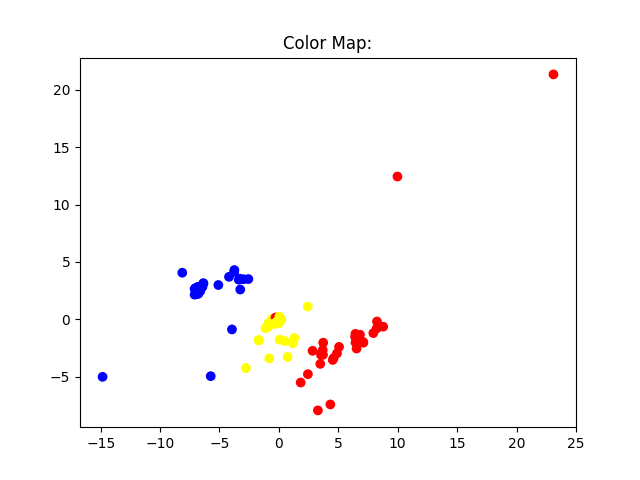}
        \includegraphics[width = 4.3cm]{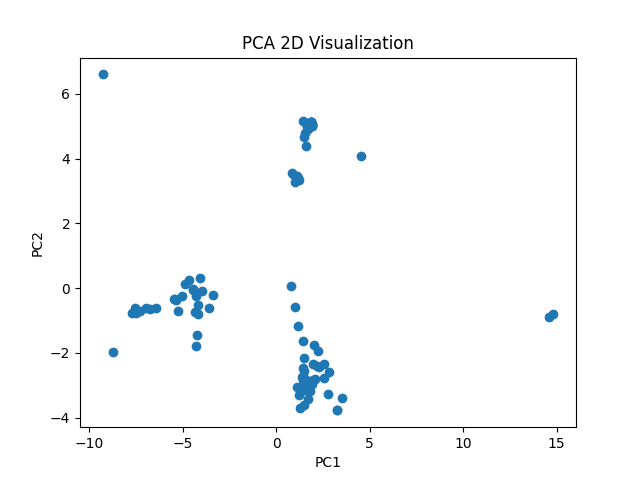}\\
        \includegraphics[width = 4.3cm]{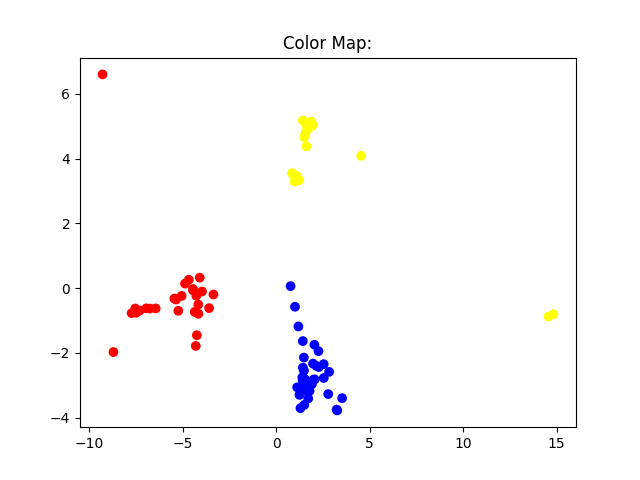}
        \includegraphics[width = 4.3cm]{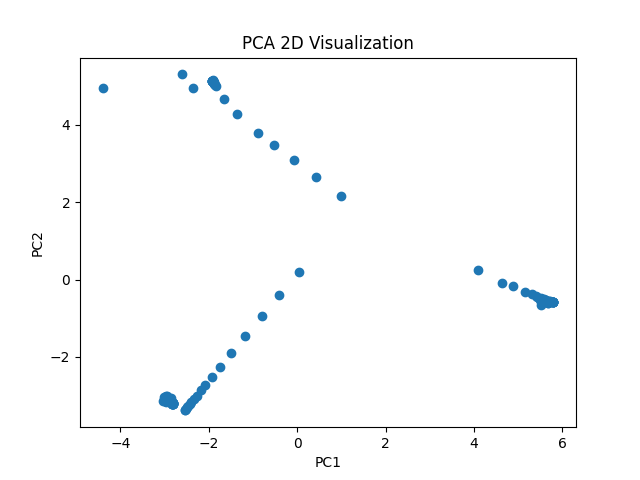}\\
    \caption{Moonset Dataset: GWA Results and PCA Visualizations (Iterations 1 to 5)}
    \label{fig:iterativemoon}
\end{figure}

Figure~\ref{fig:iterativemoon} shows the results after seven recursive iterations. 
On the initial iteration, the algorithm is able to cluster the original data into three groups. 
By using the 2D PCA results as the input for the next iteration, we are able to observe the optimal $v_i$ vectors clustering together into three distinct groups with each recursion.
This becomes clearer in the fifth iteration (bottom row) where we can clearly see the individual clusters. 
For the last two iterations, not pictured here, no significant difference is seen. 

\section{Higher Dimensional Generalization of Max $k$-Cut}
To conduct the higher dimensional generalization of our max $k$-cut algorithm, we map the original dataset with vectors of dimension $n$ into dimension $m$ where $m > n$ (recall that $n$ is the total number of points in the our original dataset).

We implement this mapping by adding rows and columns of zeroes to the original weight matrix $W$ to create a new $m$ by $m$ matrix. 
Note that, because these values are zero, no new information is added to our dataset.

Due to our results with the moonset dataset in Section~\ref{recursiveresults}, we will only be using five iterations of recursion for our results.

\subsection{Results of Higher Dimensional Generalization of GWA}
Here, we analyze the results of using both higher dimensional generalization and our recursion on the moonset dataset introduced in Section~\ref{recursiveresults}. 

\begin{figure}
\centering
        \includegraphics[width = 4.3cm]{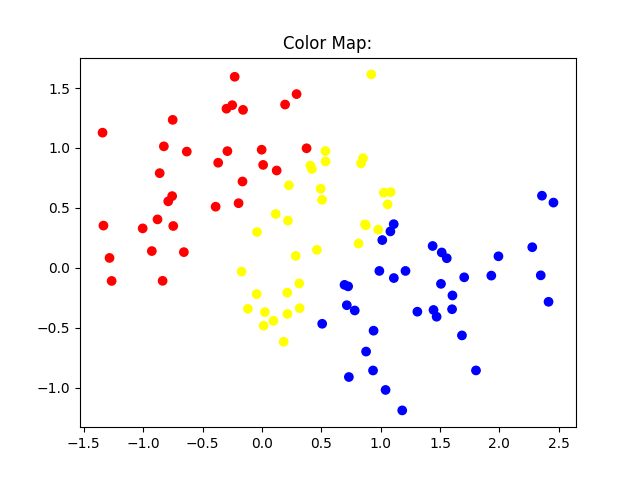}
        \includegraphics[width = 4.3cm]{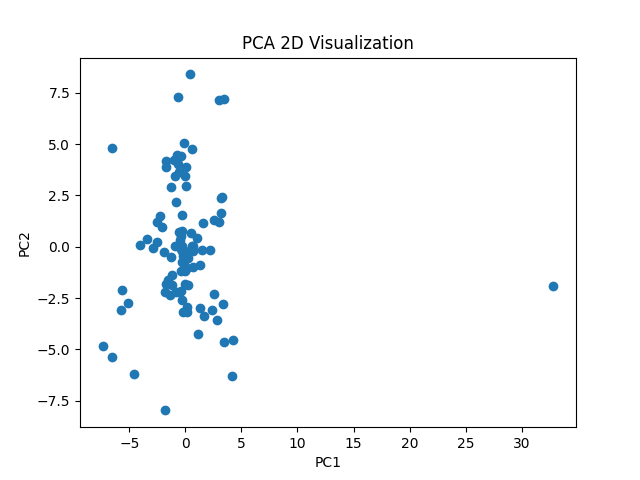}\\
        \includegraphics[width = 4.3cm]{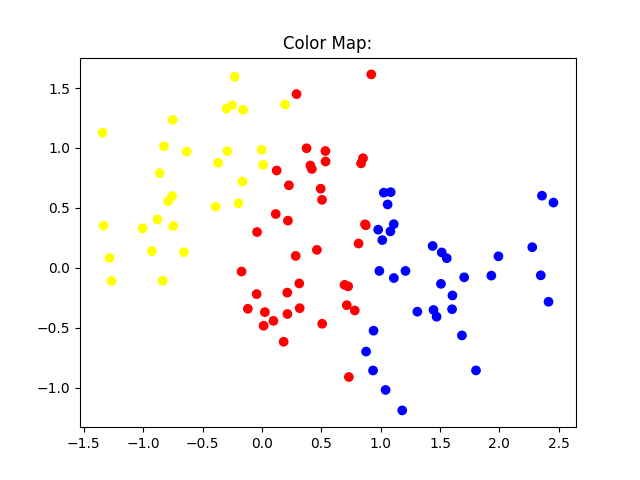}
        \includegraphics[width = 4.3cm]{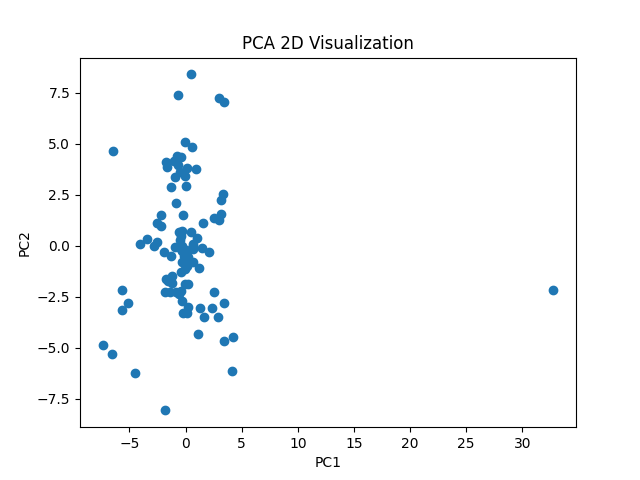}\\
        \includegraphics[width = 4.3cm]{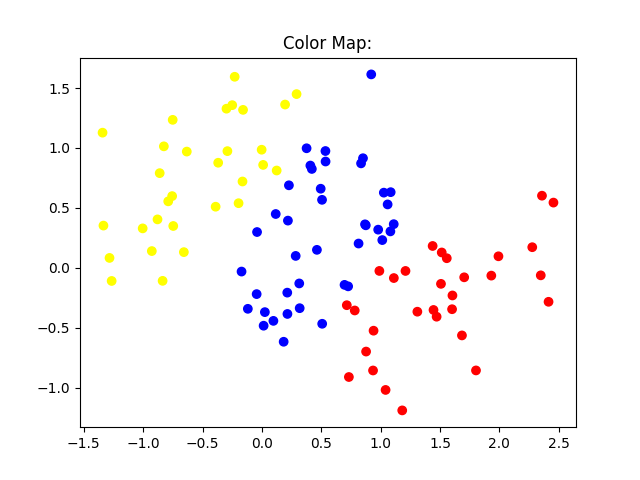}
        \includegraphics[width = 4.3cm]{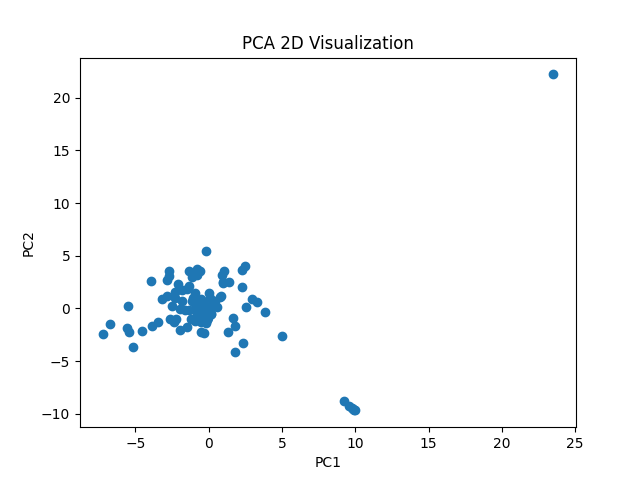}\\
        \includegraphics[width = 4.3cm]{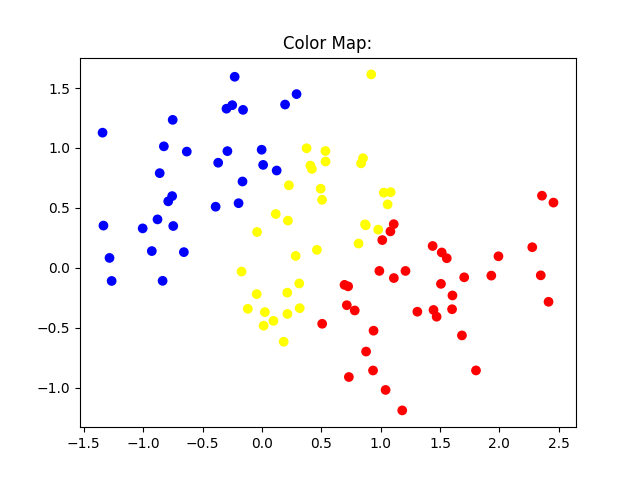}
        \includegraphics[width = 4.3cm]{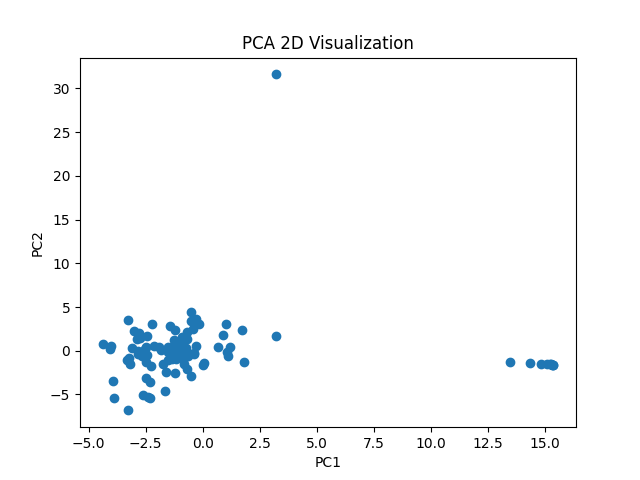}
    \caption{First Iteration: 100, 101, 106, and 109 Dimensions}
    \label{fig:recursive0}
\end{figure}

Figure~\ref{fig:recursive0} compares the results of running a single iteration of our algorithm on the moonset dataset with dimensions: 100 (the same as the size of our dataset), 101, 106, and 109. 
The few differences seen in the initial clustering of our dataset can be attributed to the stochastic nature of our algorithm as described in Steps 2 through 4 of Frieze and Jerrum's algorithm in Section~\ref{maxkcutsummary}.

Observing these iterations after applying PCA to the optimal $V$ matrix, it shows how additional dimensions can affect clustering. 
While the images for 100 and 101 dimensions (top two rows) are not significantly different, once we reach 106 dimensions (third row), the initial visualization of our $v_i$ vectors changes drastically. 
Importantly, as the dimensionality grows, the plotting scale shrinks, highlighting the tightness of the initial clusters. 
Similarly, outliers in higher dimensions are brought closer to the initial cluster with few exceptions.

\begin{figure}
\centering
        \includegraphics[width = 4.3cm]{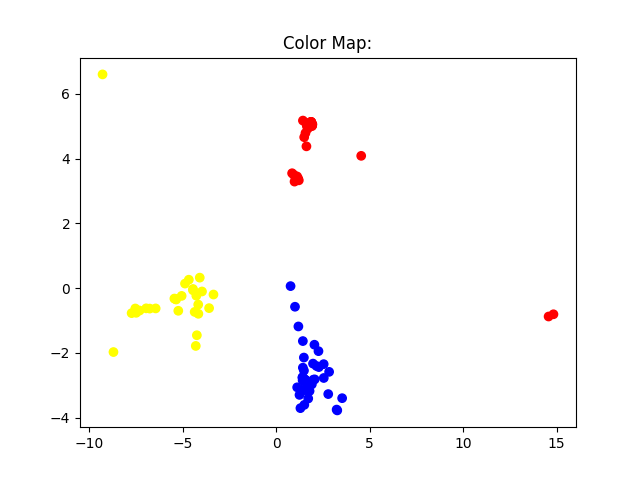}
        \includegraphics[width = 4.3cm]{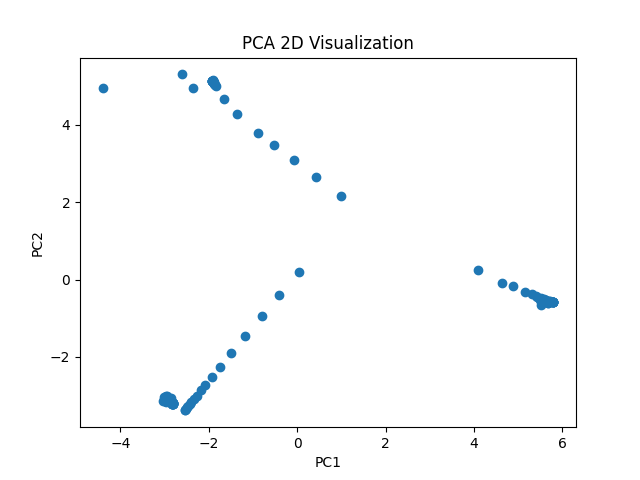}\\
        \includegraphics[width = 4.3cm]{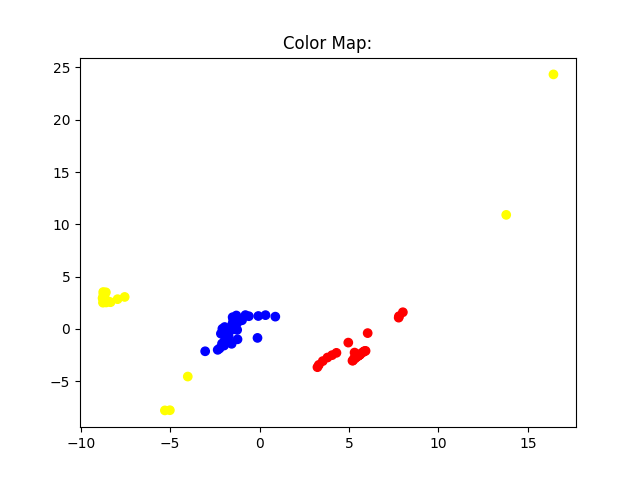}
        \includegraphics[width = 4.3cm]{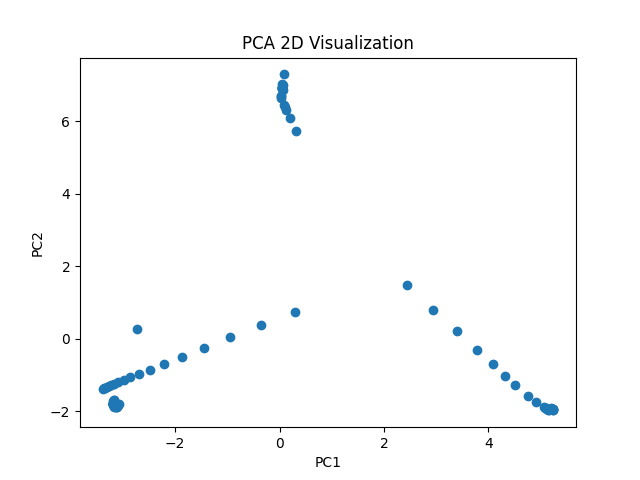}\\
        \includegraphics[width = 4.3cm]{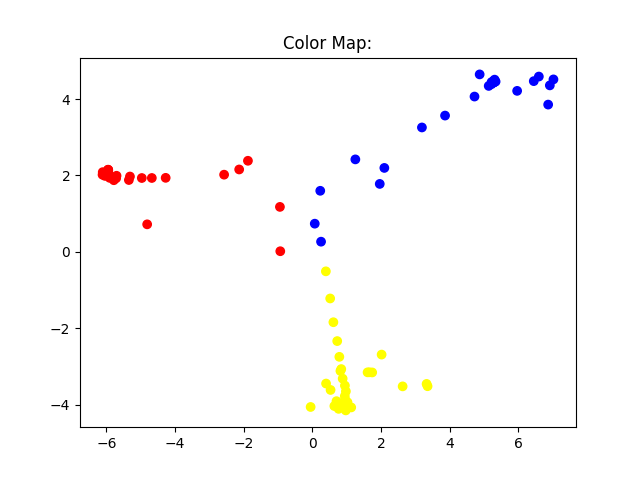}
        \includegraphics[width = 4.3cm]{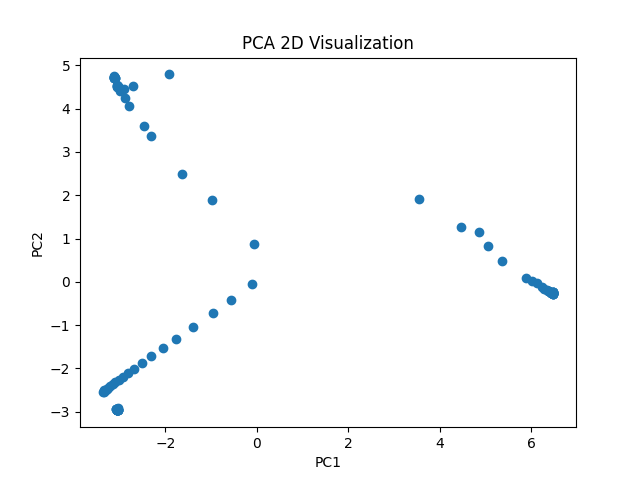}\\
        \includegraphics[width = 4.3cm]{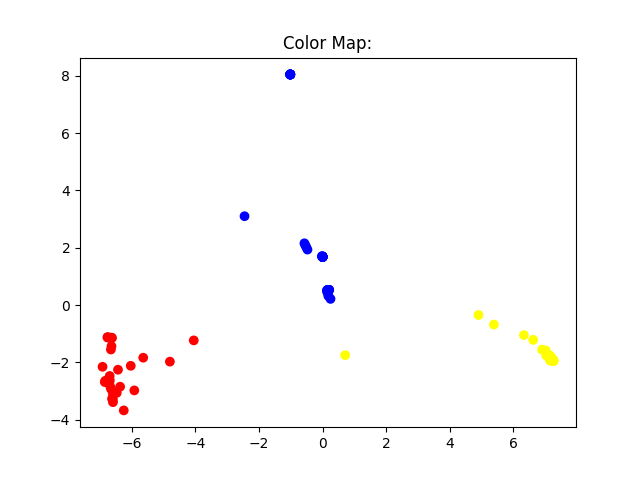}
        \includegraphics[width = 4.3cm]{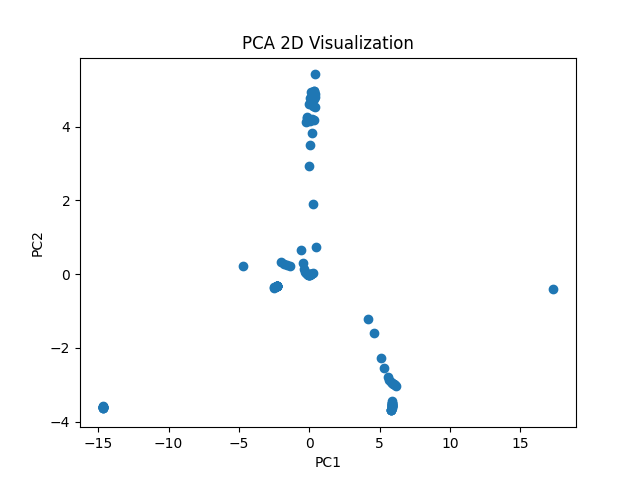}
    \caption{Fifth Iteration: 100, 101, 106, and 109 Dimensions}
    \label{fig:recursive4}
\end{figure}

We use this visualization as input for our second iteration for the recursion process. 
While there was no significant change in the first iteration, the second iteration results begins to change more drastically despite only a small dimension increase from 100 up to 109. 
Excluding all outliers across all four dimensions chosen, there is tighter grouping in higher dimensions with the clearest clusters seen in iterations 4 and 5 as seen in Figure~\ref{fig:recursive4}. 

This trend continues in subsequent iterations, resulting in Figure~\ref{fig:recursive4} which displays the results after our fifth iteration.
The PCA visualization from the fourth iteration is shown on the left column (color mapped) and the PCA visualization for the fifth iteration on the right column. 
We begin to observe clustering on iteration 4 for all dimensions as shown in the color mapped figures.
On iteration 5, these clusters become more distinct across all dimensions with the exception of 109. 
Due to the stochastic nature of this algorithm, we may be able to achieve better cluster for 109 dimensions with more recursion or clustering iterations. 
Computing limitations prevented us from this further analysis. 
Thus, we see that by the fifth iteration, the effect of dimensional relaxation results in significant differences in clustering for this dataset.

\section{Results of Clustering on Brainwave Dataset}
\label{brainwavedataset}
\begin{figure}
\centering
    \includegraphics[width = 4.3cm]{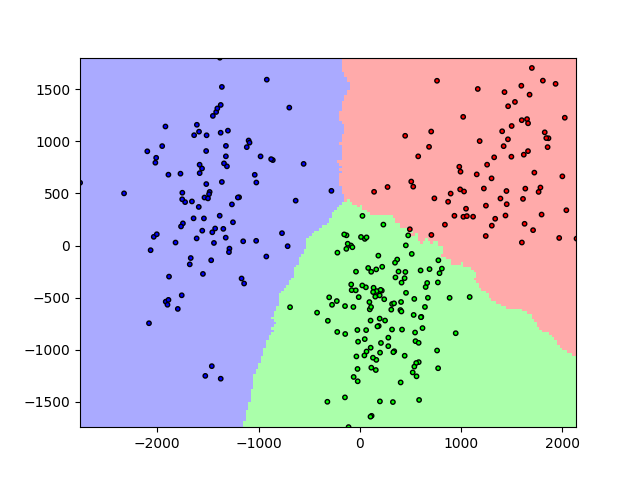}
    \includegraphics[width = 4.3cm]{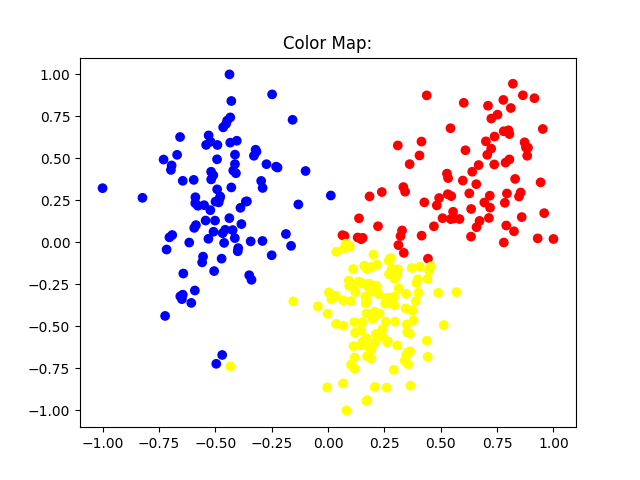}
    \caption{Brainwave Dataset: Supervised (Left) vs Unsupervised Clustering (Right)}
    \label{fig:brainwaveresults}
\end{figure}

We now apply our algorithm to a collection of datapoints detailing timestamps of spikes in brainwaves. 
While the original dataset contains 26,388 datapoints, due to the limitations of computation power and time, we test our algorithm on a set of 300 points chosen randomly from the original dataset. 
Figure~\ref{fig:brainwaveresults} compares the performance of a $k$-nearest neighbors classifier against our algorithm on this reduced dataset.
These clusters generated by both algorithms are extremely similar which is surprising given the fact that $k$-nearest neighbors is a supervised classifier and our algorithm is unsupervised. 

\begin{figure}
\centering
        \includegraphics[width = 4.3cm]{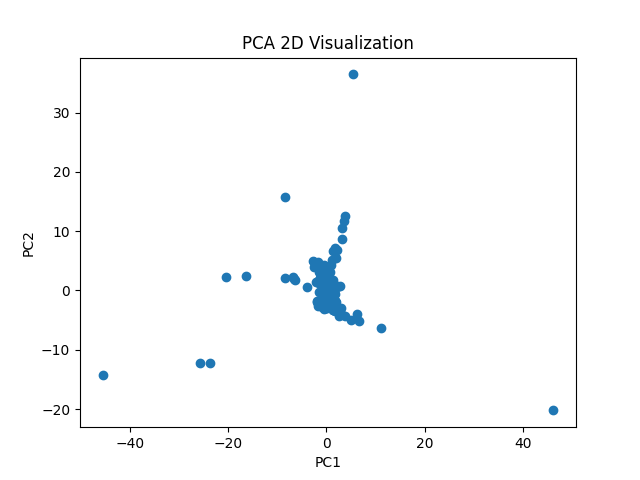}
        \includegraphics[width = 4.3cm]{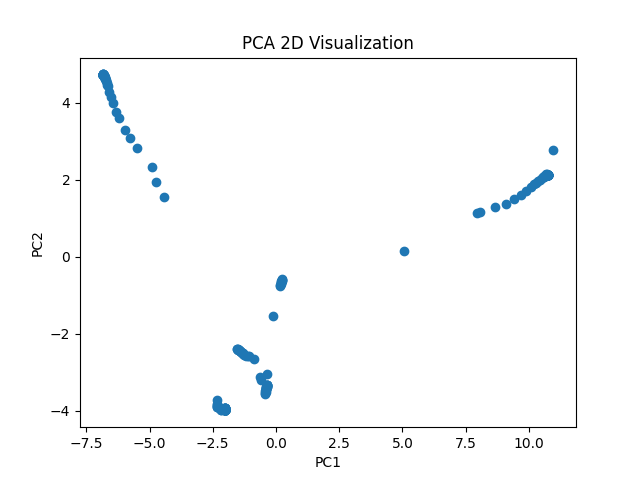}
    \caption{Brainwave Dataset: PCA Visualizations Comparison Between First (Left) and Seventh (Right) Iteration}
    \label{fig:brainwaverecursion}
\end{figure}

Figure~\ref{fig:brainwaverecursion} compares the 2D PCA visualization results of running our algorithm over a single iteration and the seventh iteration. 
Due to the presence of outliers in our initial iteration, there are significantly different scales. 
Regardless, a visual comparison of the plots show the effect that recursive iterations have on clustering as the seventh iteration provides significantly better results than the initial iteration.
Interestingly, this dataset achieves clear clusters in seven iterations whereas the dataset in Section~\ref{recursiveresults} achieves similar results in five iterations. 

\section{Vectorization of Articles}
\label{vectorization}
We begin by separating each article into its subsequent paragraphs. 
Paragraphs are identified by splitting the text when two new-line characters (“/n/n”) are seen consecutively. 
Paragraphs are then vectorized through the following procedure.

We begin by identifying a list of target words such that for each given word, we compute the conditional probability of seeing the other target words within a set window size. 
As an initial approach, we define the following target list $L$: $L$ = [amodiaquine, human, side-effect]. 
Note that we first replace any medical side effect seen within the text, i.e. ”headache”, with the word ”side-effect” before computing conditional probabilities. 
A similar procedure is used for human context words. 
The list of side effects is taken from \cite{sideeffectsebl}. 
These pre-processing steps allow for accurate probability vector calculations.

With this list, we set an arbitrary window of size $n$ centered around the target word. 
Within this window, we compute the conditional probabilities of seeing the other two target words. 
We repeat this procedure for every paragraph and can therefore use these vectors to split paragraphs of interest from irrelevant paragraphs. 
We treat the window size $n$ as a hyper-parameter and, after testing, use $n$ = 10 in the ”vectorization” process for this article specifically.

\section{Clustering Results of GWA on Amodiaquine Dataset}
\label{amodiaquinedataset}

\begin{figure}
\centering
        \includegraphics[width = 4.35cm]{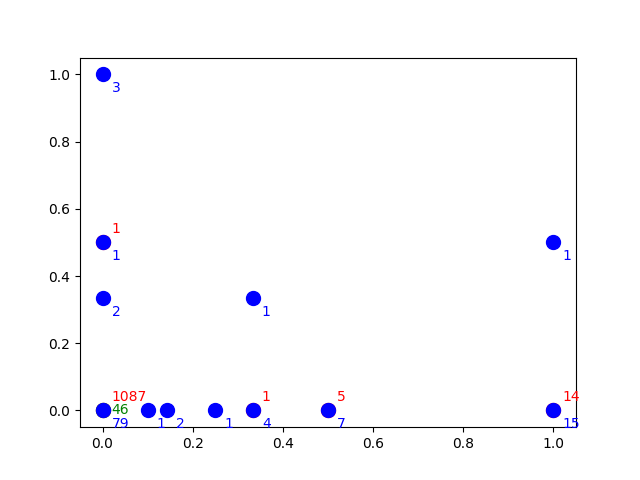}
        \includegraphics[width = 4.35cm]{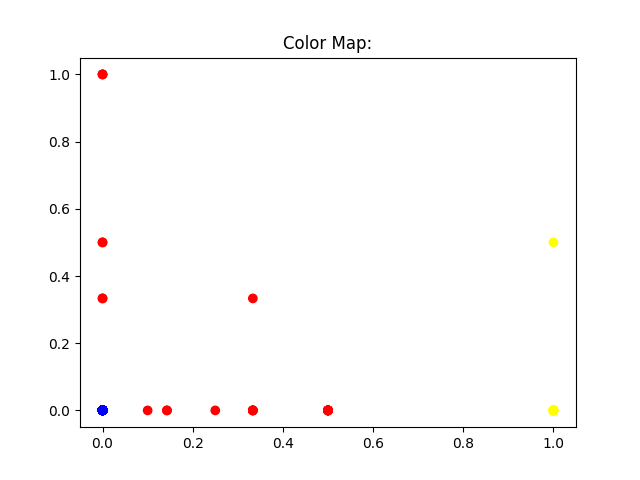}
    \caption{True Labels as Humanly-Classified (Left) and Clustering Algorithm Results (Right) with a window size equal to 10 (5 on each side)}
    \label{fig:amodiaquine}
\end{figure}

Figure~\ref{fig:amodiaquine} shows the results of our initial iteration of GWA in comparison to the true labels. 
For the graph showcasing the datapoints and their true labels on the left, we denote the number of paragraphs vectorized at a specific point by the color of the numerical label on the right of each datapoint. 
Paragraphs discussing the side effects of amodiaquine in human patients have blue text, those discussing side effects of amodiaquine in animals have green text, and those that do neither have red text. 
Note that there is overlap among all three classes, making for a skewed dataset.
The horizontal axis denotes the probability that we find a 'human' context word in the presence of 'amodiaquine' while the vertical axis denotes the probability of finding a 'side-effect' in the presence of 'amodiaquine'. 
Note that a significant portion of our desired paragraphs (blue) are at $(0,0)$ meaning that there is no mention of a side effect or a human patient within the presence of 'amodiaquine'. 
The skewness of the dataset is a result of either the window size being too small or due to the way in which paragraphs are split.
Both of these are currently being researched and further refined.

The right image in Figure~\ref{fig:amodiaquine} shows the results of our max $k$-cut algorithm on our dataset. 
Again, due to the originally skewed data, many of the points are overlapping, but there is a relatively clear separation as to where the three clusters are being centered.

\begin{figure}
\centering
        \includegraphics[width = 4.35cm]{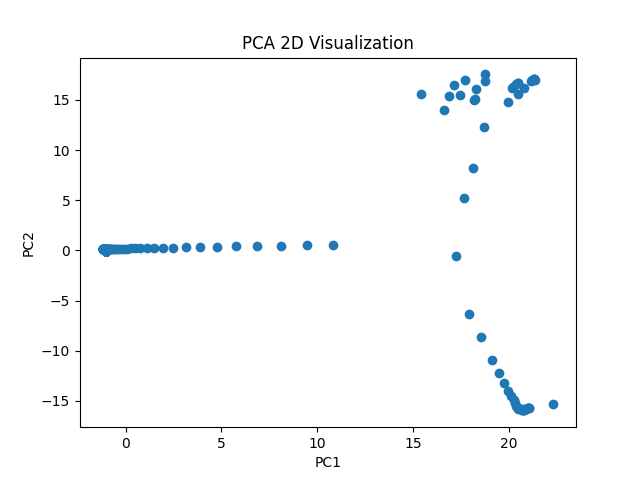}
        \includegraphics[width = 4.35cm]{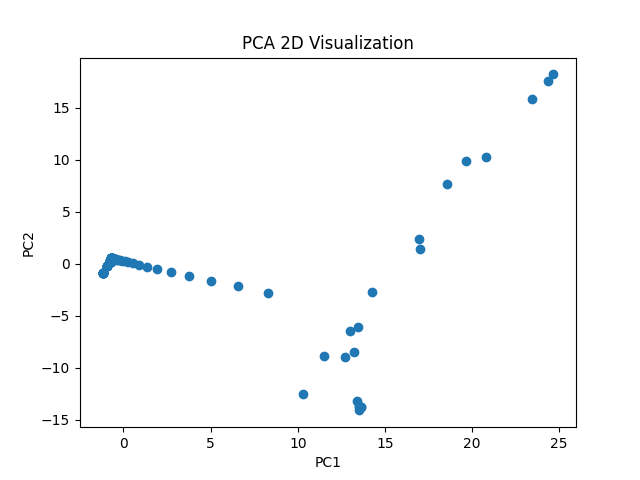}
    \caption{2D Principal Component Analysis on Full Set on 1271 (Left) and Reduced Set of 500 (Right)}
    \label{fig:vpcaanalysis}
\end{figure}

Due to computing limitations, we reduce the dataset by randomly sampling 500 points. 
After reducing the number of points, our 2D PCA visualization of our optimal $v_i$ vectors changes on the first iteration as shown in Figure~\ref{fig:vpcaanalysis}. 
Running one iteration of the algorithm on the full set shows three clear clusters as shown in the left image of the figure. 
Running one iteration of the algorithm on the reduced set shows similar results albeit with less well-defined clusters. 
This hints at the theory that, with more datapoints and computational power available to run the optimization problem, the $v_i$ vectors calculated will cluster more clearly and with fewer iterations. 

\section{Conclusions}
In conclusion, our article presents an analysis of an adaptation of a Max $k$-Cut algorithm based on work by Goemans and Williamson, additional constraints from Frieze and Jerrum, and the addition of our recursion iterations.
The effects of recursion and higher-dimensional generalization can be observed with these preliminary results for a group of three clusters, often showing better results in higher dimensions. 
In the field, this technique could be further improved to better tackle complex datasets for improved clustering in an unsupervised environment in comparison to spending a significant amount of time labeling training data for supervised clustering algorithms, especially for cases when an ample amount of data is not readily available or easily labeled for training.

Further research and refinements regarding the number of clusters and the algorithm itself are in progress.

\printbibliography
\end{document}